\documentclass[12pt, english]{article}
\usepackage[T2A]{fontenc}
\usepackage[cp1251]{inputenc}
\usepackage{babel}
\usepackage[centertags]{amsmath}
\usepackage{amsfonts}
\usepackage{amssymb}
\usepackage{amsthm}
\usepackage{newlfont}
\usepackage[active]{srcltx}

\newtheorem{theorem}{\bf Theorem}
\newtheorem{lemma}{\bf Lemma}

\begin{document}

\title{Bounds of the remainder in a combinatorial central limit theorem }

\author{Andrei N. Frolov
\\ Dept. of Mathematics and Mechanics
\\ St.~Petersburg State University
\\ St. Petersburg, Russia
\\ E-mail address: Andrei.Frolov@pobox.spbu.ru}

\maketitle

{\abstract{ We derive new bounds of the remainder in
a combinatorial central limit theorem without assumptions
on independence and existence of moments of summands.
For independent random variables our theorems imply
Esseen and Berry-Esseen type inequalities, some other new bounds and
a combinatorial central limit theorem in the case of
infinite variations.
 }


}

\section{Introduction}

Let $\|X_{ij}\|$ be a $n \times n$ matrix of independent random variables
and $\vec{\pi}=(\pi(1),\pi(2),\ldots, \pi(n))$ be a random permutation of
$1,2,\ldots, n$, independent with $X_{ij}$. Assume that $\vec{\pi}$
has the uniform distribution on the set all such permutations.
Denote
$$ S_n = \sum\limits_{i=1}^n X_{i\pi(i)}. $$

First results  on asympotical normality of $S_n$
were obtained for $P(X_{ij}=c_{ij})=1$,
$1\leqslant i, j \leqslant n$, in Wald and Wolfowits (1944).
They found sufficient conditions for that when $c_{ij}=a_i b_j$.
Noether (1949) proved that these conditions maybe
relaxed. Hoeffding (1951) considered general case
of $c_{ij}$ and obtained a combinatorial central limit theorem (CLT).
Further results on the combinatorial CLT were obtained by Motoo (1957)
and Kolchin and Chistyakov (1973).

Later investigations have been turned from limit theorems to non-asymptotic
results similar to Berry--Esseen and Esseen inequalities in
classical theory of summing of independent random variables.
Von Bahr (1976) and Ho and Chen (1978) derived bounds for the remainder
in a combinatorial CLT in the case of non-degenerated $X_{ij}$.
Botlthausen (1984) obtained Esseen type inequality for the remainder for
degenerated $X_{ij}$. The constant was not be specified
in the last paper. Further results of this type may be found in
Goldstein (2005) and Chen, Goldstein and Shao (2011).
They contain explicit constants in the inequalities.
For non-degenerated $X_{ij}$,  Esseen type inequalities were
stated by Neammanee and Suntornchost (2005), Neammanee and Rattanawong (2009)
and Chen and Fang (2012). 
These inequalities were obtained for $X_{ij}$ with finite third moments
by an application of Stein method.
At the same time, it is known that the Berry--Esseen and Esseen inequalities
maybe generalized to random variables without third moments.
This techniques for sums of independent random variables may be found in
Petrov (1995), for example.
Applying similar techniques, Frolov (2014) obtained Esseen type
bounds for the remainder in a combinatorial CLT for $X_{ij}$ with
finite variations without third moments.

In this paper, we obtain new bounds for the remainder in a combinatorial CLT
without moment assumptions. We also prove a general result in which
there are no independence assumptions.
In the case of independent random variables,
our new results generalize those in Frolov (2014). Moreover,
our results yield a combinatorial CLT for random variables
without second moments. In our example, the summands belong to
the domain of attraction of the normal law.

\section{Results}

Let $\|X_{ij}\|$ be a $n \times n$ matrix of random variables
and $\vec{\pi}=(\pi(1),\pi(2),\ldots, \pi(n))$ be a random permutation of
$1,2,\ldots, n$, where $n\geqslant 2$.
Note that we do not suppose the independence of random variables under
consideration.

Denote
$$ S_n = \sum\limits_{i=1}^n X_{i\pi(i)}. $$

For real $a_n$ and $b_n>0$, put
$$  \Delta_n  =\sup\limits_{x\in \mathbb{R}} \left|
P\left(\frac{S_n-a_n}{b_n} < x \right)-\Phi(x)\right|, $$
where $\Phi(x)$ is the standard normal distribution function.

Let $\|\mu_{ij}\|$ be a $n \times n$ matrix of real numbers and
$\|t_{ij}\|$ be a $n \times n$ matrix with
$0 < t_{ij} \leqslant +\infty$, where $n\geqslant 2$.
For $1\leqslant i, j \leqslant n$, put
$$\bar{X}_{ij} = (X_{ij}-\mu_{ij}) I\{ |X_{ij}-\mu_{ij}| < t_{ij}\},$$
where $I\{\cdot\}$ denotes the indicator of the event in brackets.
Denote
$$ \bar{S}_n = \sum\limits_{i=1}^n \mu_{i\pi(i)}
+ \sum\limits_{i=1}^n \bar{X}_{i\pi(i)},\quad
\bar{e}_n = E \bar{S}_n,\quad
\bar{B}_n= D \bar{S}_n = E \bar{S}_n^2 - (\bar{e}_n)^2,
 $$
and
$$  \bar{\Delta}_n = \sup\limits_{x\in \mathbb{R}}
\left| P\left(\frac{\bar{S}_n - \bar{e}_n}{\sqrt{\bar{B}_n}}
< x\right)-\Phi(x)\right|. $$

For all $i$ and $j$ put $p_{ij}=P(\pi(i)=j)$ and
\begin{equation*}
 q_{ij} =
  \begin{cases}
    P(|X_{ij}-\mu_{ij} | \geqslant t_{ij}|\pi(i)=j), & \text{if $p_{ij}>0$}, \\
    0, & \text{otherwise}.
  \end{cases}
\end{equation*}

Our first result is as follows.

\begin{theorem}\label{th1}
{\it 
The following inequality holds
\begin{equation}\label{be1}
\Delta_n \leqslant \bar{\Delta}_n + \Psi_n + \Theta_n + \Upsilon_n,
\end{equation}
where
$$\Psi_n =\sum\limits_{i,j=1}^n q_{ij} p_{ij}, \quad
\Theta_n =
\frac{|a_n -\bar{e}_n|}{\sqrt{2\pi} \sqrt{\bar{B}_n}},\quad
\Upsilon_n= \frac{1}{\sqrt{2\pi e}}
\max\left(\frac{\sqrt{\bar{B}_n}}{b_n}-1,
\frac{b_n}{\sqrt{\bar{B}_n}}-1 \right).
$$
}
\end{theorem}

This result is an analogue of Theorem 5.9 in Petrov (1995)
for sums of random variables.

We now turn to the main case when random variables $X_{ij}$ are independent
and permutation $\vec{\pi}$ is independent from summands and has
the uniform distribution.

For every $n\times n$ matrix $\| m_{ij} \|$, put
$$ m_{i.} =\frac{1}{n} \sum\limits_{j=1}^n m_{ij},\quad
m_{.j} =\frac{1}{n} \sum\limits_{i=1}^n m_{ij},\quad
m_{..} =\frac{1}{n^2} \sum\limits_{i,j=1}^n m_{ij}, \quad
m_{ij}^{\ast} = m_{ij} - m_{i.} - m_{.j} + m_{..}
$$
for all $i$ and $j$.

Denote $\bar{a}_{ij} = E \bar{X}_{ij}$ and $\bar{\sigma}^2_{ij} = D \bar{X}_{ij}$
for $1\leqslant i, j \leqslant n$. It is not difficult to check that
$$ \bar{e}_n= n (\bar{a}_{..}+\mu_{..}), \quad
\bar{B}_n= \frac{1}{n-1}\sum\limits_{i,j=1}^n
(\mu_{ij}^{\ast}+\bar{a}_{ij}^{\ast})^2
+\frac{1}{n} \sum\limits_{i,j=1}^n \bar{\sigma}_{ij}^2.
$$

Moreover, in this case,
$$ \Psi_n = \frac{1}{n} \sum\limits_{i,j=1}^n
P(|X_{ij}-\mu_{ij}| \geqslant t_{ij}),
$$
and Theorem \ref{th1} has the following form.

\begin{theorem}\label{th2}
{\it Assume that random variables $X_{ij}$ are independent
and permutation $\vec{\pi}$ is independent with $X_{ij}$. Suppose that
$\vec{\pi}$ has the uniform distribution on the set of all permutation
of $1,2, \ldots n$.

Then the following inequality holds
\begin{equation}\label{be2}
\Delta_n \leqslant \bar{\Delta}_n + 
\frac{1}{n} \sum\limits_{i,j=1}^n P(|X_{ij}-\mu_{ij}| \geqslant t_{ij})
+ \Theta_n +\Upsilon_n.
\end{equation}
}
\end{theorem}

There are no moment assumption in Theorems \ref{th1} and \ref{th2}.
We now consider the case of finite means.

Assume that $E X_{ij} = c_{ij}$ and
\begin{eqnarray}\label{10}
c_{i.} = c_{.j}=0,
\end{eqnarray}
for all $1\leqslant i, j \leqslant n$.
Note that this property of the matrix $\| E X_{ij}\|$
plays in a combinatorial CLT the same role that the centering at mean of
summands does in CLT.

Condition (\ref{10}) implies that $E S_n=0$ and, therefore, we take $a_n=0$.

In the sequel, we also put $t_{ij} = b_n$ for all $1\leqslant i, j \leqslant n$.

\begin{theorem}\label{th3}
{\it Assume that the conditions of theorem \ref{th2} are satisfied,
relation (\ref{10}) holds and
$\mu_{i.} = \mu_{.j}=0$ for all $1\leqslant i, j \leqslant n$.

Then there exists an absolute positive constant $A$ such that
\begin{equation}\label{be3}
  \Delta_n \leqslant
\frac{A}{n \bar{B}_n^{3/2}} \sum\limits_{i,j=1}^n
\left(|\mu_{ij}|^{3} + E|\bar{X}_{ij}|^3 \right)+
\frac{1}{n} \sum\limits_{i,j=1}^n P(|X_{ij}-\mu_{ij}| \geqslant b_n)
+ \Theta_n + \Upsilon_n.
\end{equation}
}
\end{theorem}

Note that we assume no moment conditions in Theorem \ref{th3} besides
existence of means.

Theorem \ref{th3} contains many known results and allows to derive new bounds
of remainder in a combinatorial CLT.

We start with the case of finite variations of random variables $X_{ij}$,
in which Theorem \ref{th3} yields the following result.

\begin{theorem}\label{th4}
{\it Assume that the conditions of Theorem \ref{th3} hold and
$D X_{ij} = \sigma_{ij}^2$. Put
$$ B_n = D S_n = \frac{1}{n-1}\sum\limits_{i,j=1}^n c_{ij}^2
+\frac{1}{n} \sum\limits_{i,j=1}^n \sigma_{ij}^2,
$$

Then there exists an absolute positive constant $A$ such that
\begin{equation}\label{be4}
  \Delta_n \leqslant A \left(C_n + \Lambda_n + L_n \right),
\end{equation}
where
$$ C_n = \frac{1}{n B_n^{3/2}} \sum\limits_{i,j=1}^n
|\mu_{ij}|^{3}, \quad
\Lambda_n = \frac{1}{n B_n} \sum\limits_{i,j=1}^n \alpha_{ij},\quad
L_n = \frac{1}{n B_n^{3/2}} \sum\limits_{i,j=1}^n \beta_{ij},
$$
$\alpha_{ij}=E (X_{ij}-\mu_{ij})^2 I\{|X_{ij}-\mu_{ij}| \geqslant
\sqrt{B_n}\}$ and
$\beta_{ij}=E\left|X_{ij}-\mu_{ij} \right|^3 I\{|X_{ij}-\mu_{ij}| < \sqrt{B_n} \}$
for $1\leqslant i, j \leqslant n$.
}
\end{theorem}

We would like to mention that constants $A$ are different in our theorems.
Of course, one can find them as function of the constant in
inequality (\ref{be3}). The last constant becomes from bounds
in a combinatorial CLT for summands with third moments. Unfortunately,
this constant is large now and, therefore, we do not give exact expressions here.


Theorem \ref{th4} is a generalization of Theorems 1 and 4 from Frolov (2014),
where the cases $\mu_{ij}=0$ and $\mu_{ij}=c_{ij}$
for all $1\leqslant i, j \leqslant n$ have been considered.
In the same way as in Frolov (2014), we arrive at the following result.

\begin{theorem}\label{th5}
{\it  Assume that the conditions of Theorem \ref{th4} hold.
Let $g(x)$ be a positive, even function such that
$g(x)$ and $x/g(x)$ are non-decreasing for $x>0$. Suppose that
$g_{ij}=E (X_{ij}-\mu_{ij})^2 g(X_{ij}-\mu_{ij})<\infty$
for $1\leqslant i, j \leqslant n$.

Then there exists an absolute positive constant $A$ such that
\begin{equation}\label{be5}
  \Delta_n 
\leqslant A \left( \frac{1}{ B_n^{3/2} n} \sum\limits_{i,j=1}^n |\mu_{ij}|^3 +
\frac{1}{B_n g(\sqrt{B_n}) n} \sum\limits_{i,j=1}^n g_{ij} 
\right).
\end{equation}
}
\end{theorem}

Theorem \ref{th5} includes as partial cases Theorems 2 and 5 from Frolov (2014),
where $\mu_{ij}=0$ and $\mu_{ij}=c_{ij}$ for all $1\leqslant i, j \leqslant n$,
correspondingly. For $g(x)=|x|^{2+\delta}$, $\delta \in (0,1]$,
we get the following result from Theorem \ref{th5}.

\begin{theorem}\label{th6}
{\it  Assume that the conditions of Theorem \ref{th4} hold.

Then there exists an absolute positive constant $A$ such that
$$ \Delta_n
\leqslant A \left( \frac{1}{ B_n^{3/2} n} \sum\limits_{i,j=1}^n |\mu_{ij}|^3 +
\frac{1}{B_n^{1+\delta/2} n}
\sum\limits_{i,j=1}^n E |X_{ij}-\mu_{ij}|^{2+\delta}\right),
$$
where $\delta \in (0,1]$.
}
\end{theorem}

Theorem \ref{th6} improves Theorems 3 and 6 from Frolov (2014),
where $\mu_{ij}=0$ and $\mu_{ij}=c_{ij}$ for all $1\leqslant i, j \leqslant n$,
correspondingly.

Note that Theorems \ref{th4} and \ref{th6} imply a combinatorial CLT
under Lyapunov and Lindeberg type conditions, correspondingly.

Theorems \ref{th5} and \ref{th6} may be applied to $-X_{ij}$ as well.
Nevertheless, one can derive further results from Theorem \ref{th4}
under non-symmetric conditions on distributions of $X_{ij}$ by a method
from Frolov (2014). Making use of this method, one can obtain bounds
in terms of sums of $E|X_{ij}|^{2+\delta_{ij}}$ or some
other moments depending on $i$ and $j$.

Let us turn to the case of infinite variations. In this case, Theorem \ref{th3}
also gives new results.

It is clear that we would like to put $b_n = \sqrt{\bar{B}_n}$ in this case.
The problem is that $\bar{B}_n$ depends on $b_n$. Then
consider the relation $b_n = \sqrt{\bar{B}_n}$ as an equation to
determine $b_n$. Let us show how it works on an example.

Assume that $X_{ij}$ have the same distribution with the density
$$ p(x) =
  \begin{cases}
    |x|^{-3}, & \text{if $|x|>1$}, \\
    0, & \text{otherwise}.
  \end{cases}
$$
Then $c_{ij} = \bar{a}_{ij} = 0$,
$$ \bar{\sigma}_{ij}^2 = \int\limits_{|x|<b_n} x^2 p(x) dx
= 2 \log b_n
$$
for all $1\leqslant i, j \leqslant n$ and
$$ \bar{B}_n = \frac{1}{n} \sum\limits_{i,j=1}^n \bar{\sigma}_{ij}^2 =
2 n \log b_n.
$$
It follows that the equation $b_n = \sqrt{\bar{B}_n}$ turns to
$$ b_n = \sqrt{2 n \log b_n}.
$$
It is not difficult to check that
$$ b_n \sim \sqrt{n \log n} \quad \mbox{as} \quad n \rightarrow \infty.
$$

We have
$$ P(|X_{ij}| \geqslant b_n)
= \int\limits_{|x| \geqslant b_n}  p(x) dx =
\frac{1}{b_n^2} \sim \frac{1}{n \log n}
\quad \mbox{as} \quad n \rightarrow \infty.
$$

Moreover,
$$ E |\bar{X}_{ij}|^3 = \int\limits_{|x|<b_n} x^3 p(x) dx =
2( b_n -1).
$$

Relations $b_n = \sqrt{\bar{B}_n}$ and $a_{..}=0$
imply that $\Upsilon_n=0$ and $\Theta_n=0$, correspondingly.

It follows from (\ref{be3}) that
$$ \Delta_n \leqslant A
\frac{n (b_n -1)}{ b_n^3}   +  \frac{ n}{b_n^2} =
O\left(\frac{1}{\log n}\right)
\quad \mbox{as} \quad n \rightarrow \infty.
$$

It yields that Theorem \ref{th3} gives a combinatorial CLT
with a bound for a rate of convergence. Moreover,
norming $\sqrt{n \log n}$ is determined in a similar
way as for distributions from the domain of attraction of the standard
normal law in usual CLT.

We now state a variant of a combinatorial CLT that follows from
Theorem \ref{th3}. We consider the case $\mu_{ij} = c_{ij}$.

\begin{theorem}\label{th7}
{\it Let $\{\|X_{nij}\|; 1\leqslant i, j \leqslant n, n=2,3,\ldots\}$
be a sequence of $n \times n$ matrix of independent random variables
with $E X_{nij}=c_{nij}$
and $\vec{\pi}_n=(\pi(1),\pi(2),\ldots, \pi(n))$ be random permutations of
$1,2,\ldots, n$, independent with $X_{nij}$. Assume that $\vec{\pi}_n$
has the uniform distribution on the set all permutations
of $1,2,\ldots, n$ for $n=2,3,\ldots$
Denote
$$ S_n = \sum\limits_{i=1}^n X_{n i \pi_n(i)}. $$

Assume that $c_{ni.}=c_{n.j} = 0$ for all $i$, $j$ and $n$.

Let $\{b_n\}$ be a sequence of positive constants. Put
$\bar{X}_{nij} = (X_{nij} - c_{nij}) I\{|X_{nij} - c_{nij}| < b_n \}$,
$ \bar{a}_{nij} = E \bar{X}_{nij}$, $\bar{\sigma}_{nij}^2 = D \bar{X}_{nij}$
for all $i$, $j$ and $n$.
Denote
$$ \bar{B}_n = \frac{1}{n-1} \sum\limits_{i,j=1}^n
(c_{nij}+\bar{a}_{nij}-\bar{a}_{ni.}-\bar{a}_{n.j}+\bar{a}_{n..})^2
+ \frac{1}{n} \sum\limits_{i,j=1}^n \bar{\sigma}_{nij}^2.
$$

Assume that the following conditions hold:
\begin{eqnarray*}
    &&
1)\quad \frac{1}{ b_n^{3} n} \sum\limits_{i,j=1}^n |c_{nij}|^3 \rightarrow 0
\quad \mbox{as} \quad n \rightarrow \infty,
\\ &&
2) \quad \frac{1}{n} \sum\limits_{i,j=1}^n
P(|X_{nij} - c_{nij}| \geqslant \varepsilon b_n ) \rightarrow 0
\quad \mbox{as} \quad n \rightarrow \infty \quad \mbox{for every fixed}\quad \varepsilon>0,
\\ &&
3)\quad \frac{\bar{B}_n}{b_n^2} \rightarrow 1
\quad \mbox{as} \quad n \rightarrow \infty,
\\ &&
4)\quad \frac{1}{ b_n n} \sum\limits_{i,j=1}^n |\bar{a}_{nij}| \rightarrow 0
\quad \mbox{as} \quad n \rightarrow \infty.
 \end{eqnarray*}

Then
$$ \sup_{x\in\mathbb{R}} \left|P\left(\frac{S_n}{b_n}<x\right) - \Phi(x)\right|
\rightarrow 0
\quad \mbox{as} \quad n \rightarrow \infty.
$$
}
\end{theorem}


\section{Proofs}

{\bf Proof of Theorem \ref{th1}.}
Put $p_n=b_n/\sqrt{\bar{B}_n}$, $q_n=(a_n - \bar{e}_n)/\sqrt{\bar{B}_n}$,
 \begin{eqnarray*}
    &&
\Delta_{n1} = \sup\limits_{x\in\mathbb{R}} \left|P\left(\frac{S_n-a_n}{b_n} < x \right)-
 P\left(\frac{\bar{S}_n  -a_n}{b_n} < x
\right)\right|,
   \\ &&
\Delta_{n2} = \sup\limits_{x\in\mathbb{R}} \left| P\left(\frac{\bar{S}_n - \bar{e}_n}{\sqrt{\bar{B}_n}}
< p_n x +q_n \right) -\Phi\left( p_n x +q_n \right)\right|,
    \\ &&
\Delta_{n3} = \sup\limits_{x\in\mathbb{R}} \left| \Phi\left( p_n x +q_n \right)-\Phi(x)\right|.
 \end{eqnarray*}
We have
\begin{equation*}
  \Delta_n \leqslant \Delta_{n1}+ \Delta_{n2} + \Delta_{n3}.
\end{equation*}

It is clear that $\Delta_{n2}=\bar{\Delta}_n$ and, therefore,
we will estimate $\Delta_{n1}$ and $\Delta_{n2}$.

Since
$$ S_n = \bar{S}_n + \sum\limits_{i=1}^n \left( X_{i\pi(i)}-\mu_{i\pi(i)}\right)
I\left\{|X_{i\pi(i)}-\mu_{i\pi(i)}| \geqslant t_{i\pi(i)} \right\},
$$
we have
$$ \left\{S_n < x \right\} \subset
\left\{\bar{S}_n < x \right\} \cup
\bigcup\limits_{i=1}^n
\left\{|X_{i\pi(i)}-\mu_{i\pi(i)}| \geqslant t_{i\pi(i)} \right\}.
$$
It follows that
\begin{eqnarray*}
&&
P(S_n < x) \leqslant P(\bar{S}_n < x)+ \sum\limits_{i=1}^n
P\left(|X_{i\pi(i)}-\mu_{i\pi(i)}| \geqslant t_{i\pi(i)} \right)
\\ &&
= P(\bar{S}_n < x)+\sum\limits_{i=1}^n \sum\limits_{j=1}^n
P\left(|X_{i\pi(i)}-\mu_{i\pi(i)}| \geqslant t_{i\pi(i)}, \pi(i) =j \right)
\\ &&
= P(\bar{S}_n < x)+\sum\limits_{i=1}^n \sum\limits_{j=1}^n
p_{ij} q_{ij} = P(\bar{S}_n < x)+ \Psi_n.
\end{eqnarray*}
From the other hand
$$ \left\{\bar{S}_n < x \right\} \subset
\left\{S_n < x \right\} \cup
\bigcup\limits_{i=1}^n
\left\{|X_{i\pi(i)}-\mu_{i\pi(i)}| \geqslant t_{i\pi(i)} \right\},
$$
which yields that
$$ P(\bar{S}_n < x) \leqslant P(S_n < x) + \Psi_n.
$$
It follows that
\begin{eqnarray*}
 \Delta_{n1} \leqslant  \Psi_n.
\end{eqnarray*}

The following result is a corollary of Lemma 5.2 in Petrov (1995).

\begin{lemma}\label{l1}
{\it For every real $p>0$ and $q$ the following inequality holds
\begin{equation*}
\sup\limits_{x\in\mathbb{R}}  |\Phi\left( p x +q \right)-\Phi(x)| \leqslant
\frac{|q|}{\sqrt{2 \pi}} + \frac{1}{\sqrt{2 \pi e}}
\max\left(p-1, \frac{1}{p}-1\right).
\end{equation*}
}
\end{lemma}

By Lemma \ref{l1} we get $\Delta_{n3} \leqslant \Theta_n+\Upsilon_n$.
This finishes the proof.
$\Box$

\medskip
{\bf Proof of Theorem \ref{th3}.}
We need the following known results (see, for example,  Chen and Fang (2012)).

\medskip
{\bf Theorem A.}
{\it Let $\|Y_{ij}\|$ be $n \times n$ matrix of independent random
variables such that $E Y_{ij} = \nu_{ij}$, $D Y_{ij} = \upsilon_{ij}$
and $E |Y_{ij}|^3<\infty$ for all $i$ and $j$.
Let $\vec{\pi} =(\pi(1), \pi(2), \ldots, \pi(n))$
be a random permutation of $1,2,\ldots, n$
with uniform distribution on the set of all permutations. Assume that
$\vec{\pi}$ and random variables $Y_{ij}$ are independent.

Then there exist an absolute constant $A$ such that
$$ \sup\limits_{x\in\mathbb{R}} \left|P\left(\frac{V_n - n \nu_{..}}{\sigma}<x\right)
-\Phi(x)\right| \leqslant \frac{A}{n \sigma^{3/2}}
\sum\limits_{i,j=1}^n E|Y_{ij} - \nu_{i.} - \nu_{.j} + \nu_{..}|^3,
$$
where
$$ V_n = \sum\limits_{i=1}^n Y_{i\pi(i)}, \quad
\sigma^2 = D V_n= \frac{1}{n-1} \sum\limits_{i,j=1}^n (\nu_{ij}^\ast)^2 +
\frac{1}{n-1} \sum\limits_{i,j=1}^n \upsilon_{ij}.
$$
}
\medskip

By Theorem A with $Y_{ij}= \mu_{ij}+\bar{X}_{ij}$, we have  that
$$ \bar{\Delta}_n \leqslant \frac{A}{n \bar{B}_n^{3/2}}
\sum\limits_{i,j=1}^n E\left|\mu_{ij}+\bar{X}_{ij}
-\bar{a}_{i.}-\bar{a}_{.j}+\bar{a}_{..} \right|^3.
$$
By the H\"{o}lder inequality we, get
$$ \bar{\Delta}_n \leqslant  \frac{25 A}{n \bar{B}_n^{3/2}}
\sum\limits_{i,j=1}^n \left(|\mu_{ij}|^3+E|\bar{X}_{ij}|^3
+|\bar{a}_{i.}|^3+|\bar{a}_{.j}|^3+|\bar{a}_{..}|^3\right).
$$

Making use of the Lyapunov inequality, we obtain that
$|\bar{a}_{ij}|\leqslant  (E|\bar{X}_{ij}|^3)^{1/3}$ for all $i$ and $j$.
Applying again the H\"{o}lder inequality, we write
$$ |\bar{a}_{i.}|^3 = \frac{1}{n^3}
\left|\sum\limits_{j=1}^n \bar{a}_{ij}\right|^3
\leqslant \frac{1}{n^3}
\left(\sum\limits_{j=1}^n |\bar{a}_{ij}|\right)^3
\leqslant \frac{1}{n}
\sum\limits_{j=1}^n |\bar{a}_{ij}|^3
\leqslant \frac{1}{n}
\sum\limits_{j=1}^n E |\bar{X}_{ij}|^3.
$$
It follows that
$$ \sum\limits_{i,j=1}^n |\bar{a}_{i.}|^3 \leqslant
\sum\limits_{i,j=1}^n E|\bar{X}_{ij}|^3.
$$
In the same way, we arrive at
$$ \sum\limits_{i,j=1}^n |\bar{a}_{.j}|^3 \leqslant
\sum\limits_{i,j=1}^n E |\bar{X}_{ij}|^3.
$$
Further, an application of the H\"{o}lder inequality yields that
$$ |\bar{a}_{..}|^3 = \frac{1}{n^6}
\left|\sum\limits_{i,j=1}^n \bar{a}_{ij}\right|^3
\leqslant \frac{1}{n^6}
\left(\sum\limits_{i,j=1}^n |\bar{a}_{ij}|\right)^3
\leqslant \frac{1}{n^2}
\sum\limits_{i,j=1}^n |\bar{a}_{ij}|^3
\leqslant \frac{1}{n^2} \sum\limits_{i,j=1}^n E |\bar{X}_{ij}|^3.
$$
The latter inequality implies that
$$ \sum\limits_{i,j=1}^n |\bar{a}_{..}|^3 \leqslant
\sum\limits_{i,j=1}^n E|\bar{X}_{ij}|^3.
$$

It follows that
$$ \bar{\Delta}_n \leqslant  \frac{25 A}{n \bar{B}_n^{3/2}}
\sum\limits_{i,j=1}^n \left(|\mu_{ij}|^3+ 4 E|\bar{X}_{ij}|^3\right).
$$
This bound and inequality (\ref{be2}) yield (\ref{be3}) and Theorems is proved.
$\Box$

\begin{lemma}\label{l2}
{\it Assume that the conditions of Theorem \ref{th4} hold.
Then there exists an absolute constant $A'$  such that
$$\left|1-\frac{\sqrt{\bar{B}_n}}{\sqrt{B_n}}\right|
\leqslant A' (C_n+\Lambda_n).
$$
}
\end{lemma}

{\bf Proof of Lemma \ref{l2}.}
Assume that $B_n=1$. Then
$$ \bar{X}_{ij} = (X_{ij} - \mu_{ij}) I\{ |X_{ij} - \mu_{ij}| < 1\}.
$$
Put
$$ \hat{X}_{ij} = (X_{ij} - \mu_{ij}) I\{ |X_{ij} - \mu_{ij}| \geqslant 1\}.
$$

We have
\begin{eqnarray*}
1-\bar{B}_n= B_n-\bar{B}_n
= \frac{1}{n } \sum\limits_{i,j=1}^n ( \sigma_{ij}^2 - \bar{\sigma}_{ij}^2) +
\frac{1}{n-1 } \sum\limits_{i,j=1}^n
(c_{ij}^2-
(\mu_{ij}+\bar{a}_{ij}-\bar{a}_{i.}-\bar{a}_{.j}+\bar{a}_{..})^2).
\end{eqnarray*}

Note that for all $i$ and $j$,
\begin{equation}\label{cij}
  c_{ij} = \mu_{ij}+\bar{a}_{ij} + E \hat{X}_{ij}.
\end{equation}

It follows that
\begin{eqnarray*}
\sigma_{ij}^2 - \bar{\sigma}_{ij}^2 = E(X_{ij}- \mu_{ij})^2 -
(c_{ij} - \mu_{ij})^2 - E \bar{X}_{ij}^2 + \bar{a}_{ij}^2
= E \hat{X}_{ij}^2 - 2 \bar{a}_{ij} E \hat{X}_{ij}
- (E \hat{X}_{ij})^2,
\end{eqnarray*}
for all $i$ and $j$.
Taking into account that $|\bar{a}_{ij}|<1$, we have
$$ |\sigma_{ij}^2 - \bar{\sigma}_{ij}^2| \leqslant 4 E \hat{X}_{ij}^2.
$$
Then
\begin{equation}\label{N1}
  \frac{1}{n } \sum\limits_{i,j=1}^n ( \sigma_{ij}^2 - \bar{\sigma}_{ij}^2)
  \leqslant 4 \Lambda_n.
\end{equation}

Further, applying of (\ref{cij}) implies that for all $i$ and $j$,
\begin{eqnarray*}
&&
c_{ij}^2-
(\mu_{ij}+\bar{a}_{ij}-\bar{a}_{i.}-\bar{a}_{.j}+\bar{a}_{..})^2
\\ &&
= 2 (\mu_{ij}+\bar{a}_{ij}) E \hat{X}_{ij} + (E \hat{X}_{ij})^2
+ 2 (\mu_{ij}+\bar{a}_{ij})(\bar{a}_{i.}+\bar{a}_{.j}-\bar{a}_{..})
- (\bar{a}_{i.}+\bar{a}_{.j}-\bar{a}_{..})^2.
\end{eqnarray*}

Note that
$$ \sum\limits_{i,j=1}^n \mu_{ij} \bar{a}_{i.} =
n \sum\limits_{i=1}^n \mu_{i.} \bar{a}_{i.} = 0,
$$
and, similarly,
$$ \sum\limits_{i,j=1}^n \mu_{ij} \bar{a}_{.j} = 0, \quad
\sum\limits_{i,j=1}^n \mu_{ij} \bar{a}_{..} =0.
$$
It follows that
\begin{eqnarray*}
\\ &&
\sum\limits_{i,j=1}^n
(c_{ij}^2-
(\mu_{ij}+\bar{a}_{ij}-\bar{a}_{i.}-\bar{a}_{.j}+\bar{a}_{..})^2)
\\ &&
= \sum\limits_{i,j=1}^n \left(
2 (\mu_{ij}+\bar{a}_{ij}) E \hat{X}_{ij} + (E \hat{X}_{ij})^2
+ 2 \bar{a}_{ij} (\bar{a}_{i.}+\bar{a}_{.j}-\bar{a}_{..})
- (\bar{a}_{i.}+\bar{a}_{.j}-\bar{a}_{..})^2 \right).
\end{eqnarray*}

Making use of the H\"{o}lder inequality, we have
\begin{equation}\label{N2}
  (\bar{a}_{i.}+\bar{a}_{.j}-\bar{a}_{..})^2 \leqslant 3
(\bar{a}_{i.}^2+\bar{a}_{.j}^2+\bar{a}_{..}^2).
\end{equation}

We write
\begin{eqnarray}
&& \nonumber
\sum\limits_{i,j=1}^n \bar{a}_{i.}^2 = n \sum\limits_{i=1}^n \bar{a}_{i.}^2
= \frac{1}{n} \sum\limits_{i=1}^n \left(
\sum\limits_{j=1}^n ( c_{ij} -  \mu_{ij} - E \hat{X}_{ij}) \right)^2
\\ &&  \nonumber
= \frac{1}{n} \sum\limits_{i=1}^n \left(
n c_{i.} - n \mu_{i.} - \sum\limits_{j=1}^n E \hat{X}_{ij}\right)^2
=
\frac{1}{n} \sum\limits_{i=1}^n \left(
\sum\limits_{j=1}^n E \hat{X}_{ij} \right)^2
\\ && \label{N3}
\leqslant
\sum\limits_{i,j=1}^n (E \hat{X}_{ij})^2
\leqslant
\sum\limits_{i,j=1}^n E \hat{X}_{ij}^2 = n \Lambda_n.
\end{eqnarray}
We obtain in the same way that
\begin{equation}\label{N4}
  \sum\limits_{i,j=1}^n \bar{a}_{.j}^2 \leqslant n \Lambda_n,
  \quad
\sum\limits_{i,j=1}^n \bar{a}_{..}^2 \leqslant n \Lambda_n.
\end{equation}

Further, we get by (\ref{N3}) that
\begin{equation}\label{N5}
  \sum\limits_{i,j=1}^n \bar{a}_{ij} \bar{a}_{i.} =
n \sum\limits_{i=1}^n  \bar{a}_{i.}^2 \leqslant
n \Lambda_n.
\end{equation}
It follows from (\ref{N4}) in the same way that
\begin{equation}\label{N6}
  \sum\limits_{i,j=1}^n \bar{a}_{ij} \bar{a}_{.j}
\leqslant n \Lambda_n,\quad
\sum\limits_{i,j=1}^n \bar{a}_{ij} \bar{a}_{..}
\leqslant n \Lambda_n.
\end{equation}
%

Taking into account that
$ x y \leqslant x^{3}/3 + 2 y^{3/2}/3$ for all non-negative $x$ and $y$,
we get
$$ |\mu_{ij} E \hat{X}_{ij}| \leqslant
\frac{1}{3} |\mu_{ij}|^3 +
\frac{2}{3} |E \hat{X}_{ij}|^{3/2}
\leqslant
\frac{1}{3} |\mu_{ij}|^3+
\frac{2}{3} E \hat{X}_{ij}^{2}
$$
for all $i$ and $j$. Hence
\begin{equation}\label{N7}
\sum\limits_{i,j=1}^n |\mu_{ij} E \hat{X}_{ij}| \leqslant
\frac{n C_n}{3} + \frac{2 n \Lambda_n}{3}.
\end{equation}

Since $|\bar{a}_{ij}|<1$ for all $i$ and $j$, we conclude that
\begin{equation}\label{N8}
  \sum\limits_{i,j=1}^n |\bar{a}_{ij} E \hat{X}_{ij}| \leqslant
\sum\limits_{i,j=1}^n E \hat{X}_{ij}^{2} = n \Lambda_n.
\end{equation}

It follows from (\ref{N2})--(\ref{N8}) that
$$ \frac{1}{n-1 } \sum\limits_{i,j=1}^n
(c_{ij}^2-
(\mu_{ij}+\bar{a}_{ij}-\bar{a}_{i.}-\bar{a}_{.j}+\bar{a}_{..})^2)
\leqslant \frac{2 n}{3(n-1)} C_n +
\left(\frac{4 n}{3(n-1)}+ \frac{18 n}{(n-1)} 
\right) \Lambda_n.
$$

The last inequality and (\ref{N1}) yield that
$$ |B_n-\bar{B}_n| = |1-\bar{B}_n|
\leqslant \frac{2 n}{3(n-1)} C_n +
\left(\frac{4 n}{3(n-1)}+ \frac{18 n}{(n-1)} + 4 
\right) \Lambda_n \leqslant C_n + 43 \Lambda_n,
$$
and Lemma \ref{l2} is proved for $B_n=1$.
If $B_n \neq 1$, then we apply the latter inequality to $X_{ij}/\sqrt{B_n}$,
$c_{ij}/\sqrt{B_n}$ and $\mu_{ij}/\sqrt{B_n}$.
$\Box$

\medskip
{\bf Proof of Theorem \ref{th4}.} Assume that $B_n=1$.

If $\sqrt{\bar{B}_n} \leqslant 1/2$, then by Lemma \ref{l2}
$$ \Delta_n \leqslant 1 \leqslant 2 |1-\sqrt{\bar{B}_n}|
\leqslant 2 A' (C_n+\Lambda_n).
$$
It yields (\ref{be4}) for $B_n=1$ in this case.

Assume now that $\sqrt{\bar{B}_n} \geqslant 1/2$. If
$\sqrt{\bar{B}_n} > 1$, then we have by Lemma \ref{l2} that
$$ \Upsilon_n = \frac{1}{\sqrt{2 \pi e}}
(\sqrt{\bar{B}_n}-1)
\leqslant  \frac{A'}{\sqrt{2 \pi e}} (C_n+\Lambda_n).
$$
For $\sqrt{\bar{B}_n} < 1$,  we get by Lemma \ref{l2} that
$$ \Upsilon_n = \frac{1}{\sqrt{2 \pi e}}
\left(\frac{1}{\sqrt{\bar{B}_n}}-1\right)
\leqslant \frac{1}{\sqrt{2 \pi e}} 2
(1-\sqrt{\bar{B}_n})
\leqslant \frac{2 A'}{\sqrt{2 \pi e}} (C_n+\Lambda_n).
$$

Note that
$$ \bar{X}_{ij} = (X_{ij} - \mu_{ij}) I\{ |X_{ij} - \mu_{ij}| < 1\}.
$$
Put again
$$ \hat{X}_{ij} = (X_{ij} - \mu_{ij}) I\{ |X_{ij} - \mu_{ij}| \geqslant 1\}.
$$

It is clear that
$$ \frac{1}{n} \sum\limits_{i,j=1}^n P(|X_{ij} - \mu_{ij}| \geqslant 1) =
\frac{1}{n} \sum\limits_{i,j=1}^n E I\{ |X_{ij} - \mu_{ij}| \geqslant 1\}
\leqslant \frac{1}{n} \sum\limits_{i,j=1}^n E \hat{X}_{ij}^2 = \Lambda_n.
$$
Moreover,
\begin{eqnarray*}
&&
\Theta_n = \frac{|\bar{e}_n|}{\sqrt{2 \pi} \sqrt{\bar{B}_n}}=
\frac{1}{\sqrt{2 \pi} \sqrt{\bar{B}_n} n}
\left|\sum\limits_{i,j=1}^n \bar{a}_{ij} \right|
=
\frac{1}{\sqrt{2 \pi} \sqrt{\bar{B}_n} n}
\left|\sum\limits_{i,j=1}^n E \hat{X}_{ij}
\right|
\\ &&
\leqslant
\frac{1}{\sqrt{2 \pi} \sqrt{\bar{B}_n} n}
\sum\limits_{i,j=1}^n
E \hat{X}_{ij}^2 = \frac{1}{\sqrt{2 \pi} \sqrt{\bar{B}_n} } \Lambda_n
\leqslant \frac{\Lambda_n}{\sqrt{\pi}}.
\end{eqnarray*}

These bounds imply by (\ref{be3}) that
\begin{eqnarray*}
\Delta_n \leqslant 2^{3/2} A (C_n + L_n)+
\frac{2 A'}{\sqrt{2 \pi e}} (C_n+\Lambda_n)+
\Lambda_n + \frac{\Lambda_n}{\sqrt{\pi}}.
\end{eqnarray*}
This inequality yields (\ref{be4}) for $B_n=1$.

If $B_n \neq 1$, then
we apply the result for $B_n=1$ to $X_{ij}/\sqrt{B_n}$,
$c_{ij}/\sqrt{B_n}$ and $\mu_{ij}/\sqrt{B_n}$.
$\Box$

Theorems \ref{th5} and  \ref{th6} follow from Theorem \ref{th4} in the same
way as in Frolov (2014). Details are omitted.

\medskip
{\bf Proof of Theorem \ref{th7}.}
Conditions 3) and 4) imply that $\Upsilon_n \rightarrow 0$ and
$\Theta_n \rightarrow 0$ as $n \rightarrow \infty$, correspondingly.

Take $\varepsilon>0$. We have
\begin{eqnarray*}
&&
\hspace*{-1.2\parindent}
E|\bar{X}_{nij}|^3 = E|X_{nij}-c_{nij}|^3 I\{|X_{nij}-c_{nij}|< \varepsilon b_n\}+
E|X_{nij}-c_{nij}|^3
I\{\varepsilon b_n \leqslant |X_{nij}-c_{nij}| < b_n\}
\\ &&
\hspace*{-1.2\parindent}
\leqslant \varepsilon b_n
E|X_{nij}-c_{nij}|^2 I\{|X_{nij}-c_{nij}|< b_n\}
+ b_n^3 P(|X_{nij}-c_{nij}| \geqslant \varepsilon b_n)
\\ &&
\hspace*{-1.2\parindent}
=
\varepsilon b_n (\bar{\sigma}_{nij}^2 + \bar{a}_{nij}^2)+
b_n^3 P(|X_{nij}-c_{nij}| \geqslant \varepsilon b_n)
\leqslant
\varepsilon b_n \bar{\sigma}_{nij}^2 + \varepsilon b_n^2 |\bar{a}_{nij}|+
b_n^3 P(|X_{nij} - c_{nij}| \geqslant \varepsilon b_n),
\end{eqnarray*}
for all $i$ and $j$. Hence
$$ \frac{1}{n b_n^{3}} \sum\limits_{i,j=1}^n E|\bar{X}_{nij}|^3
\leqslant \varepsilon \frac{\bar{B}_n}{b_n^2} +
\varepsilon \frac{1}{b_n n}
\sum\limits_{i,j=1}^n |\bar{a}_{nij}|
+\frac{1}{n}
\sum\limits_{i,j=1}^n P(|X_{nij} - c_{nij}| \geqslant \varepsilon b_n).
$$
This inequality and conditions 2), 3) and 4) yield that
$$ \frac{1}{n b_n^{3}} \sum\limits_{i,j=1}^n E|\bar{X}_{nij}|^3
\rightarrow 0\quad \mbox{as}\quad n \rightarrow \infty.
$$
By inequality (\ref{be3}), we arrive at the desired conclusion.
$\Box$

\bigskip
\noindent
{\bf References}
{\footnotesize

\parindent 0 mm

Bolthausen E. (1984) An estimate of the remainder in
a combinatorial central limit theorem.
Z. Wahrsch. verw. Geb. 66, 379-386.

Goldstein L. (2005) Berry-Esseen bounds for combinatorial central limit theorems
and pattern occurrences, using zero and size biasing. J. Appl. Probab. 42, 661-683.

Chen L.H.Y., Fang X. (2012) 0n the error bound in a
combinatorial central limit theorem. arXiv:1111.3159.

Chen L.H.Y., Goldstein L., Shao Q.M. (2011) Normal approximation
by Stein's method. Springer.

Ho S.T., Chen L.H.Y. (1978) An $L_p$ bounds for the remainder in
a combinatorial central limit theorem.
Ann. Probab. 6, 231-249.

Frolov A.N. (2014) Esseen type bounds of the remainder in a combinatorial CLT.
J. Statist. Planning and Inference, 149, 90-97.

Hoeffding W. (1951) A combinatorial central limit theorem. Ann. Math. Statist. 22,
558-566.

Kolchin V.F., Chistyakov V.P. (1973) On a combinatorial limit theorem.
Theor. Probab. Appl. 18, 728-739.

Motoo M. (1957) On Hoeffding's combinatorial central limit theorem.
Ann. Inst. Statist. Math. 8, 145-154.

Neammanee K., Suntornchost J. (2005) A uniform bound on
a combinatorial central limit theorem.
Stoch. Anal. Appl. 3, 559-578.

Neammanee K., Rattanawong P. (2009) A constant on a uniform bound of
a combinatorial central limit theorem. J. Math. Research 1, 91-103.

Noether G.E. (1949) On a theorem by Wald and  Wolfowitz.
 Ann. Math. Statist. 20, 455-458.

Petrov V.V. (1995) Limit theorems of probability theory. Sequences of
independent random variables. Clarendon press, Oxford.

von Bahr B. (1976) Remainder term estimate in a
combinatorial central limit theorem.
Z. Wahrsch. verw. Geb. 35, 131-139.

Wald A., Wolfowitz J. (1944) Statistical tests based on permutations of
observations. Ann. Math. Statist. 15, 358-372.

}

\end{document}